\documentclass[preprint]{svmult}

\usepackage{type1cm}        % activate if the above 3 fonts are
\usepackage{makeidx}         % allows index generation
\usepackage{graphicx}        % standard LaTeX graphics tool
\usepackage{multicol}        % used for the two-column index
\usepackage[bottom]{footmisc}% places footnotes at page bottom

\usepackage{newtxtext}       %
\usepackage{newtxmath}       % selects Times Roman as basic font
\usepackage{multirow}
\usepackage{diagbox}

\usepackage{xcolor}
\usepackage[colorlinks]{hyperref}

\definecolor{darkgreen}{rgb}{0.0,0.4,0.0}
\definecolor{darkred}{rgb}{0.6,0.0,0.0}
\definecolor{darkblue}{rgb}{0.0,0.0,0.5}
\definecolor{gray}{rgb}{0.5,0.5,0.5}
\definecolor{cyan}{rgb}{0.0,1.0,1.0}
\definecolor{darkcyan}{rgb}{0.0,0.5,0.5}
\definecolor{darkorange}{rgb}{0.8,0.4,0.0}
\definecolor{darkmargenta}{rgb}{0.5,0.0,0.5}
\definecolor{black}{rgb}{0.0,0.0,0.0}

\hypersetup{citecolor=blue, linkcolor=red, anchorcolor=darkblue,
filecolor=darkblue, urlcolor=darkblue}
\hypersetup{pdfauthor={Nicolae Suciu},pdftitle={Suciu}}

\makeindex             % used for the subject index
\begin{document}

\title*{Computational orders of convergence of iterative methods for Richards' equation}
\author{Nicolae Suciu, Florin A. Radu, Jakob S. Stokke, Emil C\u{a}tina\c{s}, Andra Malina}
\institute{Nicolae Suciu $\cdot$ Emil C\u{a}tina\c{s} $\cdot$ Andra Malina \at Tiberiu Popoviciu Institute of Numerical Analysis, Romanian Academy, Cluj-Napoca, Romania, \email{nsuciu@ictp.acad.ro, ecatinas@ictp.acad.ro, andra.malina@ictp.acad.ro}
\and
Florin A. Radu $\cdot$ Jakob S. Stokke \at Center for Modeling of Coupled Subsurface Dynamics, University of Bergen, Bergen, Norway, \email{florin.radu@uib.no, jakob.stokke@uib.no}}

%\author[2]{Florin A. Radu}
%\author[2]{Jakob S. Stokke}
%\author[1]{Emil C\u{a}tina\c{s}}
%\author[1]{Andra Malina}
%\affil[1]{Tiberiu Popoviciu Institute of Numerical Analysis, Romanian Academy, Romania}
%\affil[2]{Center for Modeling of Coupled Subsurface Dynamics, University of Bergen, Norway}
%\affil[ ]{\{nsuciu, ecatinas, andra.malina\}@ictp.acad.ro; \{jakob.stokke, florin.radu\}@uib.no}

\maketitle

\abstract*{Numerical solutions for flows in partially saturated porous media pose challenges related to the non-linearity and elliptic-parabolic degeneracy of the governing Richards' equation. Iterative methods are therefore required to manage the complexity of the flow problem. Norms of successive corrections in the iterative procedure form sequences of positive numbers. Definitions of computational orders of convergence and theoretical results for abstract convergent sequences can thus be used to evaluate and compare different iterative methods. We analyze in this frame Newton's and $L$-scheme methods for an implicit finite element method (FEM) and the $L$-scheme for an explicit finite difference method (FDM). We also investigate the effect of the Anderson Acceleration (AA) on both the implicit and the explicit $L$-schemes. Considering a two-dimensional test problem, we found that the AA halves the number of iterations and renders the convergence of the FEM scheme two times faster. As for the FDM approach, AA does not reduce the number of iterations and even increases the computational effort. Instead, being explicit, the FDM $L$-scheme without AA is faster and as accurate as the FEM $L$-scheme with AA.}

\abstract{Numerical solutions for flows in partially saturated porous media pose challenges related to the non-linearity and elliptic-parabolic degeneracy of the governing Richards' equation. Iterative methods are therefore required to manage the complexity of the flow problem. Norms of successive corrections in the iterative procedure form sequences of positive numbers. Definitions of computational orders of convergence and theoretical results for abstract convergent sequences can thus be used to evaluate and compare different iterative methods. We analyze in this frame Newton's and $L$-scheme methods for an implicit finite element method (FEM) and the $L$-scheme for an explicit finite difference method (FDM). We also investigate the effect of the Anderson Acceleration (AA) on both the implicit and the explicit $L$-schemes. Considering a two-dimensional test problem, we found that the AA halves the number of iterations and renders the convergence of the FEM scheme two times faster. As for the FDM approach, AA does not reduce the number of iterations and even increases the computational effort. Instead, being explicit, the FDM $L$-scheme without AA is faster and as accurate as the FEM $L$-scheme with AA.}

\section{Introduction}
\label{sec:intro}

The convergence of the iterative methods for Richards' equation can be proved under more or less restrictive smoothness and boundedness conditions that may not be met in some practical applications. Alternatively, the convergence may be assessed or verified numerically by investigating the behavior of the successive corrections. This could be mainly useful in the case of $L$-schemes, which do not require the evaluation of the derivatives and, as such, can be applied to problems with non-smooth coefficients and boundary conditions. Although the decay of the correction norms allows the estimation of the convergence order, this may not be enough to prove that a specific convergence behavior holds, mainly when trying to prove the linear convergence of a $L$-scheme. General methods developed for convergent sequences of real numbers allow a more in-depth convergence analysis and the evaluation of the relaxation strategies, such as AA, used to improve the iterative methods.

The convergence of an arbitrary sequence of real numbers $x_s\rightarrow x^* \in \mathbb{R}$ is characterized by the behavior of the successive errors $e_s=|x^* - x_s|$. The sequence $\{x_s\}$ converges with the (classical) $C$-order $p\geq 1$ if
\begin{equation}\label{eq:Qp}
\hspace{-0.3cm}(C)\qquad\lim\limits_{s\rightarrow\infty} \frac{e_{s+1}}{(e_s)^p}=Q_p \in
\begin{cases}(0,\infty)\;&p>1\\
             (0,1)\; &p=1,
\end{cases}
\end{equation}
which implies the (weaker) orders $Q$ and $R$\footnote{These are usually denoted by $Q_L$ and $R_L$ to distinguish from other types of $Q$ and $R$ orders. Since only the definitions (\ref{eq:p}) and (\ref{eq:R}) will be used in the following, we disregard the subscript $L$.} (see \cite{Potra-89,BEQ-90,Catinas-2019,Catinas2021}),
\begin{equation}\label{eq:p}
\hspace{-3.cm}(Q)\qquad\lim\limits_{s\rightarrow\infty} \frac{\ln e_{s+1}}{\ln e_s}=p,
\end{equation}
\begin{equation}\label{eq:R}
\hspace{-2.74cm}(R)\qquad\lim\limits_{s\rightarrow\infty} |\ln e_s|^{1/s}=p.
\end{equation}

If the limit $x^*$ is not known, the errors $e_s$ are replaced by the corrections $|x_{s+1} - x_s|$ and
(\ref{eq:Qp}-\ref{eq:R}) define ``computational orders of convergence'' $C'$, $Q'$, and $R'$. If $p>1$, computational and error-based orders of convergence are equivalent, denoted by using curly braces, and they are related by \cite{Catinas2021}
$$\{C, C'\}\underset{\nLeftarrow}{\Rightarrow}\{Q, Q'\}\underset{\nLeftarrow}{\Rightarrow}\{R, R'\}.$$
However, when $p=1$, the above equivalences do not hold. In \cite{BEQ-90} was given an example of a sequence with $C$-linear order but no $C'$-linear order. Also, in Sect. \ref{sec:num_fdm} below we present the situation of two sequences where the numerical estimates of (\ref{eq:Qp}) indicate $C$-sublinear and $C'$-linear orders of convergence.

In this framework, we examine the Newton's and $L$-scheme methods utilized in implicit FEM approaches \cite{List,Stokke2023} as well as an explicit FDM $L$-scheme \cite{Suciuetal2021,SuciuandRadu2022}.

\section{The setup of the flow problem}
\label{sec:Richards}

Combining the mass balance equation and Darcy's law, one obtains the Richards' equation for the pressure head $\psi(x,z,t)$ in variably saturated soils,
\begin{equation}\label{eq:Richards}
\frac{\partial \theta(\psi)}{\partial t}-\nabla\cdot\left[K(\theta(\psi)\nabla(\psi+z)\right]=f,
\end{equation}
where $z$ is the height against the gravitational direction, $\theta$ is the water content, $K$ stands for the hydraulic conductivity of the medium, and $f$ is a source/sink term.

Equation (\ref{eq:Richards}) is closed by relationships defining the water content $\theta(\psi)$ and the hydraulic conductivity $K(\theta(\psi))$ given by the van Genuchten-Mualem model
\begin{equation} \label{eq:theta}
\Theta(\psi) = \begin{cases} \left(1+(-\alpha \psi)^n\right)^{-m} &\psi < 0 \\
1 &\psi \geq 0,
\end{cases}
\end{equation}
\begin{equation} \label{eq:K}
K(\Theta(\psi)) = \begin{cases} K_{s} \Theta(\psi)^{\frac{1}{2}} \left[1-\left(1-\Theta(\psi)^\frac{1}{m}\right)^m \right]^2, &\psi < 0 \\
K_{s} &\psi \geq 0,
\end{cases}
\end{equation}
where $\Theta = (\theta - \theta_{r})/(\theta_{s} - \theta_{r})$ is the normalized water content, $\theta_{r}$ and $\theta_{s}$ are the residual and the saturated water content, $K_{s}$ is the saturated hydraulic conductivity, and $\alpha$, $n$, $m=1-1/n$ are model parameters depending on the soil type.

To illustrate the use of the computational orders in analyzing and comparing iterative methods for Richards' equation, we solve (\ref{eq:Richards}) closed by (\ref{eq:theta}) and (\ref{eq:K}), with parameters from \cite{Stokke2023}. We consider the domain $\Omega=\Omega_1\bigcup\Omega_2$, where $\Omega_1=[0,1]\times[0,1/4)$ and $\Omega_2=[0,1]\times[1/4,1]$ and the time interval $[0,T]$, where $T=0.003$. Sequences of corrections are stored at $T/3$, $2T/3$, and $T$.
	
As an initial condition, we choose the pressure head
\begin{equation*} \label{eq:initial}
\psi^{0}(x,z) =
\begin{cases}
-z +1/4 \;\;\;\; (x,z)\in\Omega_1, \\
-3 \quad\quad\quad\;\; (x,z)\in\Omega_2.
\end{cases}
\end{equation*}
On the top
boundary, we employ constant Dirichlet conditions, equal to the initial condition at all times. On the remaining
boundary no-flow conditions are used. We also consider the following source term
\begin{equation*} \label{eq:source}
f(x,z) =
\begin{cases}
0 \quad\quad\quad\quad\quad\quad\quad\quad\quad\quad\quad\quad\quad (x,z)\in\Omega_1, \\
0.006\cos\left(\frac{4}{3}\pi(z - 1)\sin(2\pi x)\right) \;\;\; (x,z)\in\Omega_2.
\end{cases}
\end{equation*}

\section{Linearization methods}
\label{sec:lin}

Newton's method and the $L$-scheme employed in the Galerkin FEM are reviewed in Sect.~\ref{sec:fem} and a recently introduced $L$-scheme for the FDM is shortly described in Sec.~\ref{sec:fdm}. Theoretical convergence results are presented in Sect.~\ref{sec:conv}.

\subsection{Implicit FEM approach}
\label{sec:fem}

The numerical scheme for fully discrete Richards' equation with backward Euler time discretization and Galerkin FEM reads:
given $\psi_{h,k-1}$ find $\psi_{h,k}$ such that
\begin{align}\label{eq:nonlinearrichards}
	\langle\theta(\psi_{h,k})-\theta(\psi_{h,k-1}),v_h\rangle+\Delta t\langle K(\theta(\psi_{h,k}))\nabla(\psi_{h,k}+z),\nabla v_h\rangle=\Delta t\langle f_k,v_h\rangle,
\end{align}
$\forall\; v_h\in V_h$, where $V_h$ is the Galerkin linear element space, $\langle\cdot,\cdot\rangle$ denotes the inner product, $\Delta t=T/\mathcal{K}$, $\mathcal{K}\in\mathbb{N}$, is the time step, $t_k=k\Delta t$, and $k\in\left\{1,...,\mathcal{K}\right\}$.

Newton's scheme linearizes both the $\theta(\cdot)$ and $K(\cdot)$ nonlinearities in Eq.~\eqref{eq:nonlinearrichards} as follows: given $\psi_{h,k-1},\psi_{h,k}^{s-1}\in V_h$ find $\psi_{h,k}^{s}\in V_h$ such that
\begin{align*}
	&\langle\theta '(\psi_{h,k}^{s-1})(\psi_{h,k}^{s}-\psi_{h,k}^{s-1}),v_h\rangle +\Delta t\langle K'(\theta(\psi_{h,k}^{s-1}))\nabla(\psi_{h,k}^{s}+z)(\psi_{h,k}^{s}-\psi_{h,k}^{s-1}),\nabla v_h\rangle \nonumber\\
	&=\Delta t\langle f_k,v_h\rangle -\langle\theta(\psi_{h,k}^{s-1})-\theta(\psi_{h,k-1}),v_h\rangle -\Delta t\langle K(\theta(\psi_{h,k}^{s-1}))\nabla(\psi_{h,k}^{s}+z),\nabla v_h\rangle.\label{eq:femNewton}
\end{align*}

The $L$-scheme only solves the non-linearity of $\theta(\cdot)$ by replacing the derivative $\theta '$ with a positive constant $L$. The linearization problem reads: given $\psi_{h,k-1},\psi_{h,k}^{s-1}\in V_h$ find $\psi_{h,k}^{s}\in V_h$ such that
\begin{align*} %\label{eq:femL}
	&L\langle\psi_{h,k}^{s}-\psi_{h,k}^{s-1},v_h\rangle +\Delta t\langle K(\theta(\psi_{h,k}^{s-1}))\nabla(\psi_{h,k}^{s}+z),\nabla v_h\rangle \nonumber\\
	&\qquad  =\Delta t\langle f_k,v_h\rangle -\langle \theta(\psi_{h,k}^{s-1})-\theta(\psi_{h,k-1}),v_h\rangle.
\end{align*}

In both schemes, the iterations start with the solution at the last time step, $\psi_{h,k}^{s=0}=\psi_{h,k-1}$ and are performed until the corrections to the vector $\boldsymbol\psi_{k}^{s}$ (with components $\psi_{h,k}^{s}$ in a basis of $V_h$) reach the accuracy $\|\boldsymbol\psi_{k}^{s}-\boldsymbol\psi_{k}^{s-1}\|\leq \varepsilon$ for an appropriate norm and some given tolerance $\varepsilon$. We notice that while Newton's method is 2nd order but only locally convergent, the $L$-scheme is a 1st order globally convergent method \cite{List}.

Different relaxation strategies can improve the efficiency of Newton's and $L$ schemes. In particular, with the AA the current approximation is given by a linear combination of previous approximations with coefficients that minimize the weighted average of the successive corrections \cite{Anderson1965}.

\subsection{Explicit FDM approach}
\label{sec:fdm}

Explicit $L$-schemes were introduced in \cite{Suciuetal2021}. For brevity and clarity, we consider the simpler case of the one-dimensional Richards' equation,
\begin{equation}\label{eq:Richards_1D}
\frac{\partial \theta(\psi)}{\partial t}-\frac{\partial}{\partial z}\left[K(\theta(\psi))\frac{\partial}{\partial z}(\psi+z)\right]=f,
\end{equation}
where $z\in[0,Z]$, $t\in[0,T]$.
The solution of (\ref{eq:Richards_1D}) at times $t=k\Delta t$, $t\leq T$, $k\in\mathbb{N}$, and positions $z=i\Delta z$, $i=1,\ldots, Z/\Delta z$ is first approximated by an implicit FDM scheme with backward discretization in time,
\begin{eqnarray}\label{eq:FD}
\theta(\psi_{i,k})-\theta(\psi_{i,k-1})&&=
\frac{\Delta t}{{\Delta z}^2}\left\{\left[K(\psi_{i+1/2,k})(\psi_{i+1,k}-\psi_{i,k})-
K(\psi_{i-1/2,k})(\psi_{i,k}-\psi_{i-1,k})\right] \right.\nonumber\\
&&\left.+\left(K(\psi_{i+1/2,k})-K(\psi_{i-1/2,k})\right)\Delta z\right\}+\Delta t f.
\end{eqnarray}
Next, we assign to $\psi_{i,k}$ the iteration count $s$, we add the stabilization term $L(\psi_{i,k}^{s+1}-\psi_{i,k}^{s})$ to the left-hand side of (\ref{eq:FD}) and perform iterations $s=1,2,\ldots $ starting with $\psi_{i,k}^{1}=\psi_{i,k-1}$ until the corrections to the solution $\boldsymbol\psi_{k}^{s}=(\psi_{i,k}^{s}, \ldots, \psi_{Z/\Delta z,k}^{s})$ fall under a given tolerance,
$\|\boldsymbol\psi_{k}^{s+1}-\boldsymbol\psi_{k}^{s}\|\leq \varepsilon$. The result is the following explicit $L$-scheme:
\begin{equation}\label{eq:FD-L}
\psi_{i,k}^{s+1} = \left[1-(r_{i+1/2,k}^{s}+r_{i-1/2,k}^{s})\right]\psi_{i,k}^{s}
+ r_{i+1/2,k}^{s}\psi_{i+1,k}^{s} + r_{i-1/2,k}^{s}\psi_{i-1,k}^{s} + f^{s},
\end{equation}
where
\begin{align*}
&r_{i\pm 1/2,k}^{s}=K(\psi_{i\pm 1/2,k}^{s})\Delta t/(L{\Delta z}^2), \quad r_{i\pm 1/2,k}^{s}\leq 1/2,\\
&f^{s}=\Delta t f/L+ \left(r_{i+1/2,k}^{s}-r_{i-1/2,k}^{s}\right)\Delta z-\left[\theta(\psi_{i,k}^{s})-\theta(\psi_{i,k-1})\right]/L.
\end{align*}
%Extensions to higher spatial dimensions and $L$-schemes %for parabolic transport equations are obtained in a %similar way \cite{Suciuetal2021,SuciuandRadu2022}.
\subsection{Theoretical convergence results}
\label{sec:conv}

The FEM approaches benefit from well-established theoretical results. For instance, the quadratic convergence of the Newton's method can be proved by using the Kirchhoff transformation and a regularization step \cite{Raduetal2006}. The linear convergence of the $L$-scheme is ensured by the monotonically increasing $\theta$, Lipschitz continuity of $\theta$ and $K$, bounded $\nabla\psi_{h,k}$, and the condition $L\geq L_{\theta}=\sup_{\psi}\{\theta'(\psi)\}$ \cite{List}.

Similar results for the explicit FDM $L$-scheme are not yet available. However, we can look for convergence conditions by investigating the evolution of the difference $e^{s}=\psi^{s}-\psi$ of the solutions of (\ref{eq:FD-L}) and (\ref{eq:FD}) at fixed position $i$ and time $k$. After denoting by $\frac{\Delta t}{{\Delta z}^2}\mathcal{D}(\psi^{s})$ the r.h.s. of (\ref{eq:FD}) and subtracting (\ref{eq:FD}) from (\ref{eq:FD-L}), we have
\begin{align*}
&L(e^{s+1}-e^{s})+\theta(\psi^{s})-\theta(\psi)=\frac{\Delta t}{{\Delta z}^2}[\mathcal{D}(\psi^{s})-\mathcal{D}(\psi)].
\end{align*}
Assuming (A1) $\theta'(\cdot)>0$ and (A2) $|\mathcal{D}(\psi^{s})-\mathcal{D}(\psi)|\leq L_{\mathcal{D}}|\psi^{s}-\psi|=L_{\mathcal{D}}|e^{s}|$ (strictly monotonic $\theta$ and Lipschitz continuous $\mathcal{D}$) one obtains
\[
\left|e^{s+1}\right|\leq \left|e^{s}\right| \left|1-\frac{\theta'(\psi)}{L}+\frac{L_{\mathcal{D}}\Delta t}{{L\Delta z}^2}\right|.
\]
Assuming the condition $\inf_{\psi}\{\theta'(\psi)\}>L_{\mathcal{D}}\Delta t/(\Delta z^2)$, the above relation becomes a contraction. This implies $|e^{s+1}| < |e^{s}|$, which in turn implies $\|\boldsymbol\psi_{k}^{s+1}-\boldsymbol\psi_k\| < \|\boldsymbol\psi_{k}^{s}-\boldsymbol\psi_k\|$. The condition above has to be compatible with the condition $K(\psi^{s})\Delta t/(L{\Delta z}^2)\leq 1/2$ which constrains the $L$-scheme (\ref{eq:FD-L}). These are necessary conditions for the $C$-linear convergence without being however sufficient conditions, i.e. they do not imply the existence of the limit (\ref{eq:Qp}).

\section{Numerical results}
\label{sec:num}

Computational orders of convergence are assessed by investigating the behavior of the sequence $\{x_s\}$ of successive corrections $x_s=\|\psi^s-\psi^{s-1}\|$, which are used to assess the accuracy of the approximation in the general case when the exact solution of the problem is not available. If the linearization method converges $x^*=\underset{s\rightarrow\infty}\lim x_s=0$ and $e_s=|x^*-x_s|=x_s$. The error-based orders $C$,  $Q$, and $R$ of the sequence $\{x_s\}$ computed according to (\ref{eq:Qp}-\ref{eq:R}), if they exist\footnote{In a multidimensional setting the $C$ orders of convergence are norm-dependent \cite{Catinas-2019}. However, the numerical examples presented below indicate the existence of the limit from the definition (\ref{eq:Qp}).}, are in fact computational orders $C'$, $Q'$, and $R'$ for the sequence $\{\psi^s\}$. To proceed, we use the definitions (\ref{eq:p}) and (\ref{eq:R}) to assess a convergence order $p$, which is further used in (\ref{eq:Qp}) to check whether the classical $C$-convergence of the sequence $\{x_s\}$ takes place. The linearization methods presented in Sect.~\ref{sec:lin} are tested in Sects.~\ref{sec:num_fem}~and~\ref{sec:num_fdm} on the benchmark problem formulated in Sec.~\ref{sec:Richards} with a fixed tolerance $\varepsilon = 10^{-7}$. In Sects.~\ref{sec:num_fdm} we also assess error-based orders of convergence for a test problem with known solution modeling variably saturated flow fully coupled with surfactant transport.

\subsection{FEM Newton's method and $L$-scheme}
\label{sec:num_fem}

\begin{figure}[h]
\begin{minipage}[c]{0.49\linewidth} %\centering
\includegraphics[width=\linewidth]{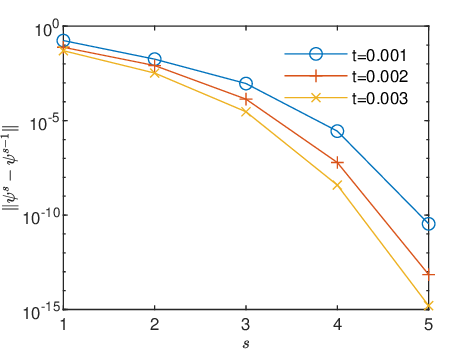}
\caption{\label{fig:fem_convN}FEM: Norms of successive corrections provided by Newton's method.}
\end{minipage}
\hspace{0.3cm}
\begin{minipage}[c]{0.49\linewidth} %\centering
\includegraphics[width=\linewidth]{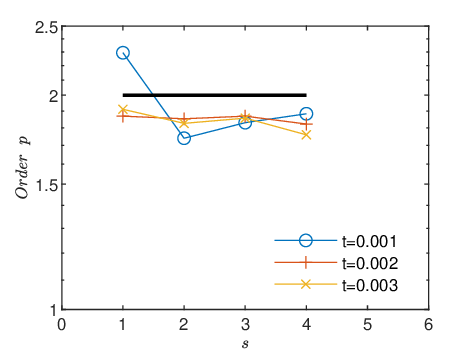}
\caption{\label{fig:fem_QordN}FEM: Estimation of $Q$-order $p$ for
the convergence of Newton's method.}
\end{minipage}
\end{figure}

Newton's method converges after only a few iterations (Fig.~\ref{fig:fem_convN}). The estimate $p(s)=\ln e_{s+1}/\ln e_s\approx 2$ shown in Fig.~\ref{fig:fem_QordN} suggests the $Q$-quadratic convergence of the sequence $\{x_s\}$, according to definition (\ref{eq:p}). The sequence of successive approximations $\{\psi^s\}$ of the Newton's method is then $Q'$-quadratic convergent.

\begin{figure}[h]
\begin{minipage}[c]{0.49\linewidth} %\centering
\includegraphics[width=\linewidth]{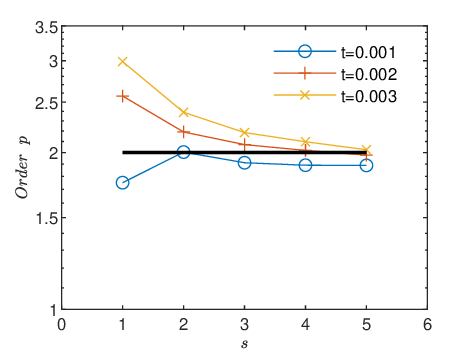}
\caption{\label{fig:fem_convR}FEM: Estimation of $R$-order 2 for Newton's method.}
\end{minipage}
\hspace{0.3cm}
\begin{minipage}[c]{0.49\linewidth} %\centering
\includegraphics[width=\linewidth]{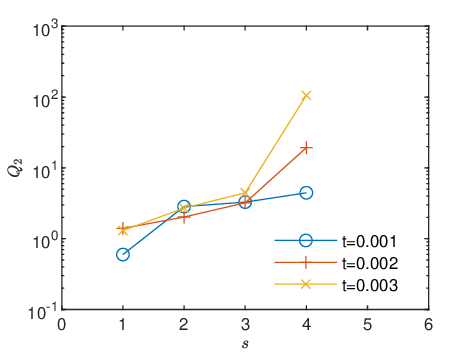}
\caption{\label{fig:fem_Q2}FEM: Estimation of the quotient $Q_2$ from definition (\ref{eq:Qp}) for Newton's method.}
\end{minipage}
\end{figure}
The estimate $p(s)=|\ln e_s|^{1/s}\approx 2$ presented in Fig.~\ref{fig:fem_convR} also indicates the $R'$-quadratic convergence according to definition (\ref{eq:R}), sustaining the hypothesis of the $Q'$-quadratic order. The finite values $Q_2(s)=e_{s+1}/e_s^2$ shown in Fig.~\ref{fig:fem_Q2} are consistent with the definition (\ref{eq:Qp}) and we may expect the $C'$-quadratic convergence. Provided that the equivalence of computational and error-based orders for $p>1$ in $ \mathbb{R}^n$ \cite{Catinas2021} holds in the normed space of the solutions, the sequence $\{\psi^s\}$ also possesses $C$, $Q$, and $R$ quadratic orders. Longer sequences necessary to validate these results could only be obtained by doing multi-precision computations, because the $e_s$ values from Fig.~\ref{fig:fem_convN} are already close to the machine epsilon in double precision.

\begin{figure}[h]
\begin{minipage}[t]{0.45\linewidth}\centering
\includegraphics[width=\linewidth]{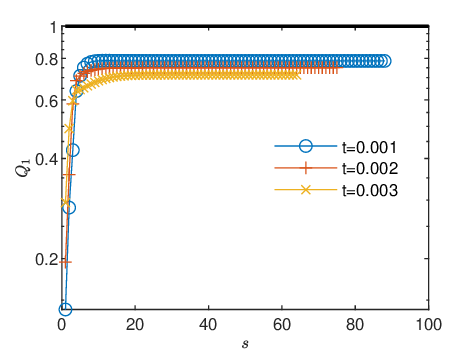}
\caption{\label{fig:fem_Q1L}FEM: Estimation of the quotient $Q_1$ from definition (\ref{eq:Qp}) for $L$-scheme ($L=0.15$).}
\end{minipage}
\hspace{0.3cm}
\begin{minipage}[t]{0.45\linewidth}\centering
\includegraphics[width=\linewidth]{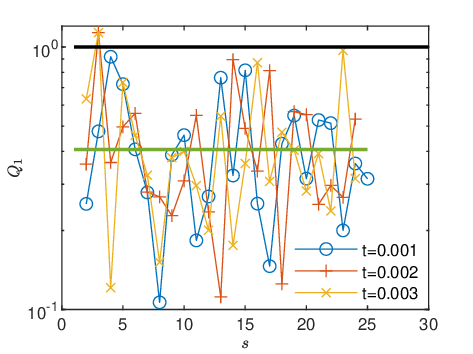}
\caption{\label{fig:fem_Q1LAA}FEM:  Estimation of the quotient $Q_1$ from definition (\ref{eq:Qp}) for $L$-scheme with AA.}
\end{minipage}
\end{figure}

Similarly, using (\ref{eq:p}) and (\ref{eq:R}) we first verified the $Q'$- and $R'$-convergence and we identified the order $p=1$ for the sequence $\{\psi^s\}$ of solutions provided by the FEM $L$-scheme. The behavior of the quotient $Q_1(s)=e_{s+1}/e_s$ shown in Fig.~\ref{fig:fem_Q1L} provides a strong indication of the $C'$-linear convergence according to the definition (\ref{eq:Qp}) with $p=1$. The scheme benefits from the use of AA which, for the setup of Sect.~\ref{sec:Richards}, halves the number of iterations and renders the convergence two times faster (Fig.~\ref{fig:fem_Q1LAA}).

\subsection{FDM $L$-scheme}
\label{sec:num_fdm}

\begin{figure}[h]
\begin{minipage}[c]{0.49\linewidth} %\centering
\includegraphics[width=\linewidth]{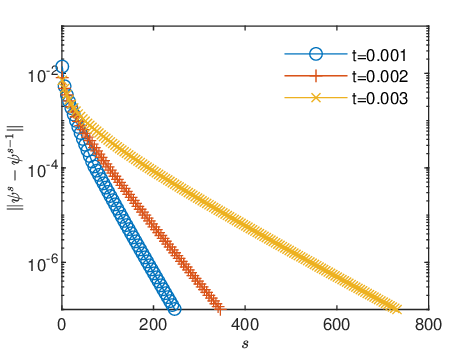}
\caption{\label{fig:fdm_convL}FDM: Norms of successive corrections provided by the $L$-scheme ($L=0.5$).}
\end{minipage}
\hspace{0.3cm}
\begin{minipage}[c]{0.49\linewidth} %\centering
\includegraphics[width=\linewidth]{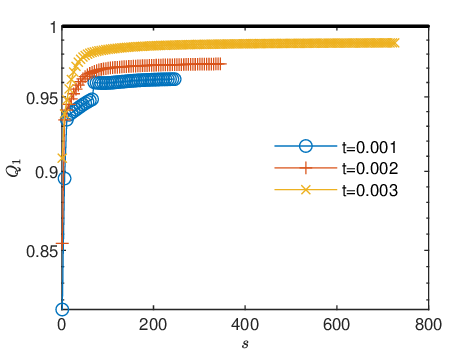}
\caption{\label{fig:fdm_Q1L}FDM: Estimation of the quotient $Q_1$ from definition (\ref{eq:Qp}) for $L$-scheme.}
\end{minipage}
\end{figure}

The number of iterations for the explicit FDM $L$-scheme (Fig.~\ref{fig:fdm_convL}) is one order of magnitude larger than for the implicit FEM $L$-scheme without AA. However, the computing time is generally one order of magnitude smaller than for implicit schemes in solving similar problems (see \cite{Suciuetal2021}). The $Q'$- and $R'$-linear convergence and the order $p=1$ are assessed in the same way as for the FEM $L$-scheme. Estimates of $Q_1$ shown in Fig.~\ref{fig:fdm_Q1L} indicate the $C'$-linear convergence. Applying AA to the explicit $L$-scheme only results in increasing the computing time, without any significant improvement of the explicit FDM $L$-scheme.

Error-based orders of convergence for linearly convergent $L$-schemes also can be inferred when the exact solution of the problem is known. As an example, we consider a test problem for fully coupled flow and surfactant transport with  the manufactured exact one-dimensional solution $\psi_{m}(z,t)=-tz(1-z)+z/4$, for the pressure, and $c_m(z,t)=tz(1-z)+1$, for the concentration, with $\theta(\psi)=1/(1-\psi-c/10)$ and $K(\theta(\psi))=\psi^2$ (the one-dimensional version of the problem considerd in
\cite[Sect. 5.2.1]{Suciuetal2021}). The problem is solved by utilising the FDM $L$-scheme for Richards' equation coupled with an explicit global random walk $L$-scheme for the advection-diffusion equation. The parameter $L=100$ is the same for both schemes and the stopping criterion uses errors instead of corrections\footnote{We note that once this stopping criterion for errors is satisfied, the corresponding corrections for both $\{\psi^s\}$ and $\{c^s\}$ are of the order $10^{-8}$.}, $\max(\|\psi_{k}^{s}-\psi_{m,k}\|,\|c_{k}^{s}-c_{m,k}\|)\leq 10^{-3}$.

\begin{figure}[h]
\begin{minipage}[t]{0.45\linewidth}\centering
\includegraphics[width=\linewidth]{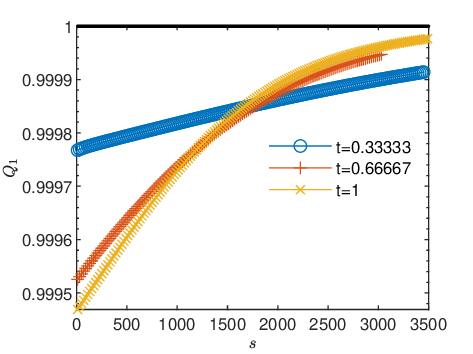}
\caption{\label{fig:fdm_Q1Le}FDM: Quotient $Q_1$ for $L$-scheme estimated with errors $\|\psi_{k}^{s}-\psi_{m,k}\|$.}
\end{minipage}
\hspace{0.3cm}
\begin{minipage}[t]{0.45\linewidth}\centering
\includegraphics[width=\linewidth]{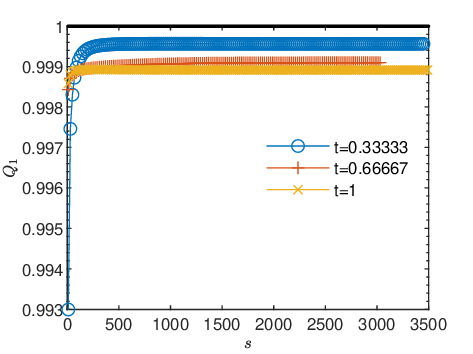}
\caption{\label{fig:fdm_Q1Lc}FDM: Quotient $Q_1$ for $L$-scheme estimated with corrections $\|\psi_{k}^{s}-\psi_{k}^{s-1}\|$.}
\end{minipage}
\end{figure}

\begin{figure}[h]
\begin{minipage}[t]{0.45\linewidth}\centering
\includegraphics[width=\linewidth]{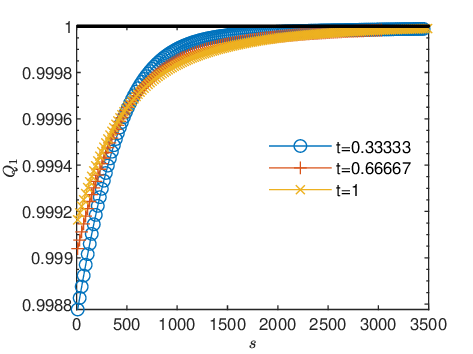}
\caption{\label{fig:fdm_Q1Le_c}FDM: Quotient $Q_1$ for $L$-scheme estimated with errors $\|c_{k}^{s}-c_{m,k}\|$.}
\end{minipage}
\hspace{0.3cm}
\begin{minipage}[t]{0.45\linewidth}\centering
\includegraphics[width=\linewidth]{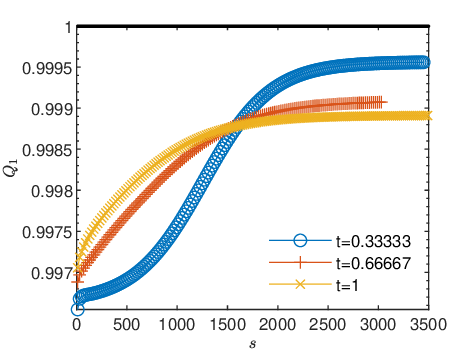}
\caption{\label{fig:fdm_Q1Lc_c}FDM: Quotient $Q_1$ for $L$-scheme estimated with corrections $\|c_{k}^{s}-c_{k}^{s-1}\|$.}
\end{minipage}
\end{figure}

The quotient $Q_1(s)$ computed with errors of the sequence
$\{\psi^s\}$ shows a tendency towards the $C$-sublinear convergence, i.e., $Q_1\rightarrow 1$ (Fig.~\ref{fig:fdm_Q1Le}). The evolution of the corresponding quotient of corrections, with $Q_1<1$, indicates the $C'$-linear convergence (Fig.~\ref{fig:fdm_Q1Lc}). The sequence of concentrations, $\{c^s\}$, behaves similarly, i.e. $C$-sublinear convergence of the errors and $C'$-linear convergence of the corrections (Figs.~\ref{fig:fdm_Q1Le_c}~and\ref{fig:fdm_Q1Lc_c}).

Thanks to the manufactured solution, we also assessed the estimated order of convergence (EOC) of the numerical FDM solution towards the exact solution $(\psi_{m},c_{m})$, by successively halving $\Delta z$, with $\max(\|\psi_{k}^{s}-\psi_{k}^{s-1}\|,\|c_{k}^{s}-c_{k}^{s-1}\|)\leq 10^{-6}$ as stopping criterion for the iterations of the coupled $L$-schemes (Table~\ref{table1}).

\begin{table}[h]
\begin{tabular}{ c c c c c c c  c c c c}
  \hline
   & $\Delta z$ & & $\|\psi-\psi_{m}\|$ &  & EOC &  & $\|c-c_{m}\|$ &  & EOC & \\
  \hline
   &  1.00e-1  & & 4.59e-02 &  & --   &  & 7.64e-03 &  & --   & \\
   &  5.00e-1  & & 1.14e-02 &  & 2.00 &  & 1.95e-03 &  & 1.97 & \\
   &  2.50e-2  & & 2.95e-03 &  & 1.95 &  & 4.38e-04 &  & 2.16 & \\
   &  1.25e-2  & & 9.97e-04 &  & 1.57 &  & 8.99e-05 &  & 2.28 & \\
  \hline
\end{tabular}
\caption{\label{table1}Orders of convergence of the numerical FDM solution.}
\end{table}

\section{Discussion}
\label{sec:disc}

The convergence of the linearization methods for Richards' equation can be investigated with methods developed for sequences of real numbers. To do that, we follow the recipe: first, use the definition of the $Q$- and $R$-order convergence to identify the convergence order $p$, then, check whether the sequence $\{x_s\}=\{\|\psi^s-\psi^{s-1}\|\}$ converges with the classical $C$-order $p$.

We have to keep in mind that the error-based convergence orders $C$, $Q$, and $R$ for the sequence of positive numbers $\{x_s\}$ correspond to the computational orders $C'$, $Q'$, and $R'$ of the sequence of solutions $\{\psi^s\}$ provided by the linearization scheme. Error-based orders for the linearization scheme only can be inferred when exact (e.g., manufactured) solutions of the flow problem are available.

The numerical analysis of the iterative methods for nonlinear systems of equations involved in modeling processes in porous media presented in the literature aim at showing the convergence to zero of the error of the scheme. Numerically inferred computational or error-based orders of convergence complete the characterization of the convergence and may be useful in several ways: to verify the theoretically predicted convergence and to check the convergence in cases where the conditions of the convergence theorems are hard to meet, to evaluate the efficiency of the relaxation strategies, or to enable comparisons of the convergence speed for different linearization procedures.


\begin{thebibliography}{99}

\bibitem{Anderson1965}Anderson, D.G.: Iterative procedures for nonlinear integral equations. JACM {\bf 12}(4), 547-560, (1965)

\bibitem{BEQ-90}
Beyer, W. A., Ebanks, B. R., Qualls, C. R.: 
Convergence rates and convergence-order profiles for sequences. Acta Appl. Math. {\bf 20}, 267--284, (1990)

\bibitem{Catinas-2019}
C\u{a}tina\c{s}, E.: A survey on the high convergence orders and computational convergence orders of sequences. Appl. Math. Comput. {\bf 343}, 1--20, (2019)

\bibitem{Catinas2021}C\u{a}tina\c{s}, E.: How many steps still left to x*? SIAM Rev. {\bf 63}(3), 585--624, (2021)

\bibitem{List}List, F., Radu, F.A.: A study on iterative methods for solving Richards' equation. Comput. Geosci. {\bf 20}(2) 341-353, (2016)

\bibitem{Potra-89}Potra, F. A.: On {Q}-order and {R}-order of convergence. J. Optim. Theory Appl. {\bf 63}(3), 415--431, (1989)

\bibitem{Raduetal2006} Radu, F.A., Pop, I.S., Knabner, P.: On the convergence of the Newton method for the mixed finite element discretization of a class of degenerate parabolic equations. In: Numerical Mathematics and Advanced Applications {\bf 42}, 1192-1200, (2006)

\bibitem{Suciuetal2021}Suciu, N., Illiano, D., Prechtel, A., Radu, F.A.: Global random walk solvers for fully coupled flow and transport in saturated/unsaturated porous media. Adv. Water Resour. {\bf 152}, 103935, (2021)

\bibitem{SuciuandRadu2022}Suciu, N., Radu, F.A.: Global random walk solvers for reactive transport and biodegradation processes in heterogeneous porous media. Adv. Water Res., {\bf 166}, 104268, (2022)

\bibitem{Stokke2023} Stokke, J.S., Mitra, K., Storvik, E., Both, J.W., Radu, F.A.: An adaptive solution strategy for Richards' equation. Comput. Math. Appl. {\bf 152}, 155-167, (2023)

\end{thebibliography}
\end{document}